\documentclass[12pt]{elsarticle}
\usepackage{graphicx,epsfig}
\usepackage{amsmath,amssymb,euscript}
\usepackage{latexsym}
\usepackage{color}
\usepackage{subfigure}
\journal{Arxiv}
\newtheorem{theorem}{\bf Theorem}[section]

\newtheorem{lemma}[theorem]{\bf Lemma}
\newtheorem{proposition}[theorem]{\bf Proposition}

\newtheorem{corollary}[theorem]{\bf Corollary}
\newtheorem{remark}[theorem]{\bf Remark}

\newcommand\al{{\alpha}}

\newcommand\be{{\beta}}
\def\NN{{\mathbb N}}

\def\PP{{\mathbb P}}
\def\C{{\mathcal{C}}}
\def\O{{\mathcal{O}}}
\def\H{{\mathcal{H}}}

\def\E{{\mathcal{E}}}

\newcommand{\Proofend}{\begin{flushright} $\square$ \end{flushright}}
\begin{document}
\begin{frontmatter}
\title{Filtered  integration rules for finite Hilbert transforms}


\author[mymainaddress,mymainaddress1]{Donatella Occorsio}
\ead{donatella.occorsio@unibas.it}

\author[mymainaddress,mymainaddress1]{Maria Grazia Russo\corref{cor1}}
\ead{mariagrazia.russo@unibas.it}

\author[mymainaddress1]{Woula Themistoclakis}
\ead{woula.themistoclakis@cnr.it}

\address[mymainaddress]{Department of Mathematics, Computer Science and Economics,\\ University of Basilicata, Via dell'Ateneo Lucano 10, 85100 Potenza, Italy }

\address[mymainaddress1]{C.N.R. National Research Council of Italy, IAC Institute for Applied Computing ``Mauro Picone'',\\ Via P. Castellino, 111, 80131 Napoli, Italy. }

\begin{abstract}
A product quadrature rule, based on the filtered de la Vall\'ee Poussin polynomial approximation, is proposed for evaluating the finite Hilbert transform in $[-1,1]$. Convergence results are stated in weighted uniform norm for functions belonging to suitable Besov type subspaces. Several numerical tests are provided, also comparing the rule with other formulas known in literature.
\end{abstract}

\end{frontmatter}
\section{Introduction}
The numerical computation of the Hilbert transform of a function plays an important role in several fields, since many mathematical models in applied sciences lead to it
(see e.g.  \cite{kalandya,king} and the references therein). Depending on the specific application, we may consider bounded or unbounded integration domains where the function may have several degrees of smoothness (see e.g. the treatise by F. King in  two volumes \cite{king} on many aspects of the Hilbert transform and its possible variants). Here we consider the case of locally continuous functions $f$ on the reference domain $[-1,1]$ and, for any Jacobi weight
\[
u(x)=v^{a,b}(x):=(1-x)^a (1+x)^b,\qquad a,b>-1,
\]
we focus on the numerical approximation of the finite Hilbert transform of $f$ (also known as the ``density function'') defined as follows
\begin{equation}\label{iniziale}
\mathcal{H}^uf(t):=
\int_{-1}^1\frac {f(x)}{x-t}u(x)dx=\lim_{\epsilon\rightarrow 0}\int_{|x-t|\ge \epsilon}\frac {f(x)}{x-t}u(x)dx,\qquad -1<t<1,
\end{equation}
being $\H\equiv \H^u$ in the case $u=1$.

Incidentally, such transform appears in Cauchy singular integral equations, which, in turn, arise  in several mathematical models (see for instance \cite{Kress,PS,Mandal}).
Moreover, $\mathcal{H}^u f$ is related to hypersingular operators that can be defined as derivatives of the Hilbert transform
and appear in integro-differential equations (see e.g. \cite{Berthold,JuLu-JCAM97,air,Monegato_siam, PS,Cauchy,Monegato_prandtl,CCJL,DBoccothem, DeBonis2021}).

Due to the relevance of the problem, there exists a wide literature on the quadrature rules for $\H^uf$. We cite for instance \cite{CMR,giova,crimas_gau,hasegawa,monegato_computing,monegato_numerishe, sinc, spline}. In particular, concerning the unweighted case, quadrature rules based on equidistant nodes of $[-1,1]$ have been recently considered in \cite{Filoccothem}.

In this paper, as system of nodes,  we shall consider the zeros of Jacobi polynomials $\{p_n(w)\}_n$ associated to suitable Jacobi weights $w=v^{\alpha,\beta}$, say
\begin{equation}\label{Xn}
X_n(w):=\{x_k:=x_{n,k}(w) \ : \ k=1,\ldots,n\},\qquad n\in\NN.
\end{equation}

Based on Jacobi systems of nodes,  several quadrature rules are known in  literature (see e.g. \cite{giova,crimas_prod,davis,monegato_computing} and the references therein). In particular, for any $n\in \NN$, we recall the product rule (L-rule) obtained  by replacing $f$ with the Lagrange polynomial $L_n(w,f)$ interpolating $f$ at  $X_n(w)$ (see e.g. \cite{monegato_computing,CMR,crimas_prod}).

Here we propose a similar product rule (that we will call VP-rule), but instead of the Lagrange polynomial, we employ the filtered de la Vall\'ee Poussin (VP) polynomial of $f$ \cite{woula_NA,them_L1}
\begin{equation*}\label{VP}
V_n^m(w,f,x)=\sum_{k=1}^n f(x_k)\Phi_{n,k}^m(x),\qquad n>m\in\NN, \qquad x\in [-1,1],
\end{equation*}
based on the same $n$ nodes required for $L_n(w,f)$. Differently from  $L_n(w,f)$,   the polynomial $V_n^m(w,f)$ does not interpolate $f$ except in some special cases \cite{woula_NA,Occo_Them, Mata2020}. Its main  characteristic  is the dependence on the additional degree--parameter $0<m<n$ which determines the action ray of the  VP filter
\vspace{.1cm}

{}\hspace{1cm}
\setlength{\unitlength}{.65cm}
\begin{picture}(7,3.5)
\put(0.3,0.3){\circle*{.125}}\put(0,0){\makebox(0,0){0}}
\put(0.3,0.3){\vector(1,0){7}}\put(7.7,0.3){\makebox(0,0){}}
\put(0.3,0.3){\vector(0,1){2.5}}\put(0,3.3){\makebox(0,0){}}
\put(0,2){\makebox(0,0){1}}\put(0.3,2){\circle*{.125}}
\put(0.3,2){\color{red}\line(1,0){2.5}}
\put(2.8,2){\color{red}\line(1,-1){1.7}}
\put(2.8,0){\makebox(0,0){n-m}}\put(2.8,0.3){\circle*{.125}}
\put(2.8,0.6){\circle*{.08}}\put(2.8,0.9){\circle*{.08}}
\put(2.8,1.2){\circle*{.08}}\put(2.8,1.5){\circle*{.08}}\put(2.8,1.8){\circle*{.08}}
\put(4.5,0){\makebox(0,0){n+m}}\put(4.5,0.3){\circle*{.125}}
\put(4.5,2.2){VP filter $\mu_{n}^m(x)=\left\{\begin{array}{ll}
1 & \mbox{if}\quad 0\le x\le n-m,\\ [.1in]
\frac{n+m-x}{2m} & \mbox{if}\quad
n-m< x< n+m.
\end{array}\right.$}
\end{picture}
\vspace{.5cm}\newline
Such filter function defines the fundamental VP polynomials
as follows \cite{ThBa}
\begin{equation*}\label{f-VP}
\Phi_{n,k}^m(x)=\lambda_k\sum_{j=0}^{n+m-1}\mu_n^m(j)p_j(w,x_k)p_j(w,x),
\end{equation*}
being
\begin{equation}\label{chris}
\lambda_k:=\lambda_{n,k}(w)= \left( \sum_{j=0}^{n-1}[p_j(w,x_k)]^2\right)^{-1},\qquad k=1,2,\ldots,n,
\end{equation}
the Christoffel numbers related to the weight $w$.

Note that in the limiting case $m=0$ the previous  VP  polynomial $V_n^m(w,f)$ coincides with the Lagrange polynomial $L_n(w,f)$ and therefore it yields the same quadrature rule.

The aim of the present paper is to show that, by using $V_n^m(w,f)$ with $0<m<n$, we can take advantage of the additional parameter $m$ which  can be suitably modulated in order to improve the quadrature error of the L-rule.

Indeed, it is already known (see e.g. \cite{woula_NA,Occo_Them_AMC,Mata2020,Occo_Them}) that corresponding to suitable choices of $0<m<n$, the polynomial $V_n^m(w,f)$  provides a pointwise approximation better than the one offered by $L_n(w,f)$, especially in presence of Gibbs phenomena. In fact, if  $f$ presents some   ``pathologies''  (peaks,   cusps, etc.)  localized in  isolated points,  the  Gibbs phenomenon affects $L_n (w,f)$ with overshoots and oscillations that spill over  the whole interval, but such phenomenon appears strongly reduced by using  $V_n^m (w,f)$.

The experimental results will show that such an improvement is inherited by  VP-rules, that for these kinds of density functions may provide a  performance better  than the L-rule.

In our experiments we will also compare VP-rules with another class of quadrature formulas based on Jacobi zeros, i. e. the Modified Gaussian rules (shortly denoted by MG rules), proposed in \cite{crimas_gau} and also studied in \cite{CMR}.
Such formulas provide   higher performance than L-rule, but we will show some cases where VP-rules produce an even better quadrature error.
Moreover, we remark that differently from both L and VP-rules,  MG rules require a variable
number of nodes depending on the singular point $t$ and consequently their employment in numerical methods for  singular integral equations appears more complicated than the previous product rules.

From a theoretical point of view, in this paper we provide new estimates of the VP polynomial approximation in some Besov type spaces characterized by a Dini--type condition. Using such estimates we prove the convergence of the proposed quadrature rule  and state several asymptotic bounds for the quadrature error.

The outline of the paper is the following. Section 2 contains preliminary notations and results concerning the approximation spaces. Moreover,  in a dedicated subsection, the mapping properties of the Hilbert transform in such spaces are studied. Section 3 is devoted to the VP approximation, while the new quadrature rule is introduced in Section 4, where we state the main results. Finally, Section 5 concerns the numerical tests while Section 6 is devoted to the proofs.

\section{Notations and preliminary results}
Consider a fixed Jacobi weight
$$v(x)=v^{\gamma,\delta}(x):=(1-x)^\gamma(1+x)^\delta, \qquad x\in (-1,1), \quad \gamma,\delta\ge 0.$$
Throughout the paper the space of all functions $f$ continuous on $(-1,1)$ and satisfying
\begin{equation}\label{lim}
\lim_{x\rightarrow +1}f(x)v(x)=0\quad \mbox{if $\gamma>0$},
\qquad \mbox{and}
\qquad
\lim_{x\rightarrow -1}f(x)v(x)=0\quad \mbox{if $\delta>0$},
\end{equation}
is denoted by $C^0_v$ and equipped with the norm
\[
\|f\|_{C^0_v}:=\|fv\|_\infty= \max_{x\in [-1,1]}|f(x)|v(x).
\]
It is well-known (see for instance \cite{mastromilobook}) that the Weierstrass approximation theorem holds in this Banach space and we have
\[
f\in C^0_v \ \Longleftrightarrow \ \lim_{n\rightarrow \infty} E_n(f)_v=0,
\]
 where $E_n(f)_v$  denotes the error of best approximation of $f\in C^0_v$ in the space $\PP_n$ of all algebraic polynomials of degree at most $n$, namely
\[
E_n(f)_v:=\inf_{P\in\PP_n}\|(f-P)\|_{C^0_v}.
\]
The rate of convergence of such error as $n\rightarrow\infty$ depends on the smoothness of the function $f$ and it is well characterized by the following moduli of smoothness introduced in \cite{DT} by Z.Ditzian and V.Totik
\[
\Omega_\varphi^k(f,t)_v=\sup_{0<h\leq t} \|v\Delta^k_{h\varphi}f\|_{L^\infty[-1+2h^2k^2, 1-2h^2k^2]}, \]
\[\omega_\varphi^k(f,t)_v=\Omega_\varphi^k(f,t)_v+\inf_{q\in\PP_{k-1}} \|(f-q)v\|_{L^\infty[-1,-1+t^2k^2] } +\inf_{q\in\PP_{k-1}} \|(f-q)v\|_{L^\infty[1-t^2k^2,1] },\]
 where
$$
\Delta^k_{h\varphi}(f,x)=\sum_{i=0}^k (-1)^i \binom{k}{i} f\left(x+\frac{kh}{2}\varphi(x)-ih\varphi(x)\right).
$$
In fact, for $n\in \mathbb{N}$ sufficiently large (say $n\geq n_0$) and $t>0$ sufficiently small (say $t\leq t_0$), the following Jackson and Stechkin type inequalities hold (see \cite{DT}, \cite{LutherRusso})
\begin{eqnarray}
\label{direct}
&&\hspace{-1.1cm}E_n(f)_v\le \C \omega_\varphi^k\left(f,\frac 1n\right)_v\le \C \int_0^\frac 1n\frac{\Omega_\varphi^k(f,t)_v}tdt, \qquad \C\ne\C(n,f),
\\
\label{converse}
&&\hspace{-1.1cm}\Omega_\varphi^k(f,t)_v\le \omega_\varphi^k(f,t)_v\le \C t^k \sum_{n=0}^{[1/t]}(n+1)^{k-1}E_n(f)_v, \qquad \C\ne\C(n,f,t),
\end{eqnarray}
where throughout the paper we use $\C$ in order to denote a positive constant, which may have different values at  different occurrences, and we write $\C\ne\C(n,f,\ldots)$ to mean that $\C>0$ is independent of $n,f, \ldots$.

Note that the previous direct and inverse results yield, for instance, the following equivalences
\begin{equation*}\label{equi_best}
E_n(f)_v=\O(n^{-r})\Longleftrightarrow \omega_\varphi^k(f,t)_v=\O(t^r)\Longleftrightarrow \Omega_\varphi^k(f,t)_v=\O(t^r)
\end{equation*}
holding for any $0<r<k\in\NN$.

In literature, several approximation spaces have been introduced in order to classify the smoothness of the functions $f$ w.r.t. the decay of $E_n(f)_v$ as $n\rightarrow \infty$. In particular, we recall the following ones
\begin{eqnarray*}
&&B_{r}(v):=\left\{f\in C^0_v \ :\sum_{n=1}^\infty (n+1)^{r-1} E_n(f)_v<\infty\right\}, \qquad r>0,
\\
&&\mbox{equipped with}
\quad \|f\|_{B_{r}(v)}:=\|fv\|_\infty+\sum_{n=1}^\infty (n+1)^{r-1} E_n(f)_v.
\end{eqnarray*}
These Banach spaces belong to the class of Besov type space $B_{r,q}^\infty(v)$, introduced in \cite{DTpapero} and  characterized by the following norm
\[
\|f\|_{B_{r,q}^\infty(v)}:=\|fv\|_\infty+\left\{
\begin{array}{ll}
\displaystyle\left(\sum_{n=1}^\infty \left[(n+1)^{r-\frac 1q} E_n(f)_v\right]^q\right)^\frac 1q,
& \mbox{if $\quad 1\le q<\infty$},\\ & \\
\displaystyle\sup_{n>0}(n+1)^r E_n(f)_v, & \mbox{if $\quad q=\infty$}.
\end{array}\right.
\]
Hence $B_{r}(v)\equiv B_{r,1}^\infty(v)$.

In the sequel, for $A,B>0$ we will write $A\sim B$ meaning that $\C^{-1}A\le B\le \C A$, with $\C>0$ independent of the main parameters in $A, B$.
By virtue of (\ref{direct})--(\ref{converse}) the following equivalence between the norms holds true \cite[Th. 3.1]{DTpapero},
\begin{eqnarray}
\label{equi-B}
\|f\|_{B_{r}(v)}&\sim & \|fv\|_\infty+\int_0^1\frac{\Omega_\varphi^k(f,t)_v}{t^{r+1}}dt,\qquad k> r>0,
\end{eqnarray}
being the constant $\C$, involved in such equivalence, only depending on $r$.
Moreover, in (\ref{equi-B})  the main part modulus $\Omega_\varphi^k$ can be also replaced by the complete modulus $\omega_\varphi^k$ or by the following K--functional
\[
K_\varphi^k(f,t^k)_v:=\inf_{g^{(k-1)}\in AC_{loc}}\left[\|(f-g)v\|_\infty + t^k\|g^{(k)}\varphi^k v\|_\infty\right]
\]
where $AC_{loc}$ denotes the class of locally absolutely continuous functions on $[-1,1]$ (i.e. absolutely continuous on any $[c,d]\subset (-1,1)$). We also recall that \cite{DT}
\[
K_\varphi^k(f,t^k)_v\sim \omega_\varphi^k(f,t)_v.
\]
For any $r>0$, the previous spaces $B_r(v)$ are compactly embedded in a larger space $B_0(v)$ that is defined \cite{LutherRusso} as the ``limit'' case of the Besov spaces $B_r(v)$ as $r\rightarrow 0$, namely
\[
B_{0}(v):=\left\{f\in C^0_v \ :\sum_{n=1}^\infty \frac{E_n(f)_v}{n+1}<\infty\right\},
\]
and it is equipped with the norm
\[
\|f\|_{B_0(v)}:=\|f v\|_\infty+\sum_{n=1}^\infty \frac{E_n(f)_v}{n+1}.
\]
Also in this case, by virtue of (\ref{direct})--(\ref{converse}), we have the following equivalences
\begin{equation}\label{normB0-equi}
\|f\|_{B_{0}(v)}\sim \|fv\|_\infty+\int_0^1\frac{\omega_\varphi^k(f,t)_v}{t}dt\sim \|fv\|_\infty+\int_0^1\frac{K_\varphi^k(f,t^k)_v}{t}dt.
\end{equation}
Moreover, the elements $f$ of $B_0(v)$ can be easily characterized with the help of the classical modulus of continuity of $g=fv$, defined as
\[
\omega(g,t):=\sup_{x,y\in [-1,1], \ |x-y|\le t}|g(x)-g(y)|.
\]
Indeed a practical criterion to check whether $f\in C^0_v$ belongs to $B_0(v)$ or not is the following \cite[Th.2.6]{Uwe2005}
\[
f\in B_0(v) \Longleftrightarrow \int_0^1\frac{\omega(fv,t)}t dt<\infty,
\]
and we also have the equivalence between the norms
\begin{equation*}\label{normB0-equi1}
\|f\|_{B_{0}(v)}\sim \|fv\|_\infty+\int_0^1\frac{\omega(fv,t)}{t}dt.
\end{equation*}
For all $r\ge0$, if $v=1$ ($\gamma=\delta=0$) we use the notation $B_r$ instead of $B_r(v)$, and have the following
\begin{lemma}\label{lem-u}
For all $r\ge 0$ and any $v=v^{\gamma,\delta}$ with $\gamma,\delta\ge 0$ we have that $f\in B_r(v)$ iff it is $(fv)\in B_r$ and (\ref{lim}) holds. Moreover we have
\begin{equation}\label{equi-lemma}
\|f\|_{B_r(v)}\sim \|fv\|_{B_r}, \qquad \forall f\in B_r(v).
\end{equation}
\end{lemma}
Finally the following result states some asymptotic bounds for the error of the best polynomial approximation in $B_r(v)$.
\begin{lemma}\label{lemmaEm}
For any $n$ sufficiently large (say $n\geq n_0$) it results
\begin{equation}\label{stimaEm}
E_n(f)_v\leq \mathcal{C} \left\{\begin{array}{ll} \displaystyle
\frac{\|f\|_{B_0(v)}}{\log n},  & \mbox{if $f\in B_0(v)$}, \\
& \\
\displaystyle \frac{\|f\|_{B_{r}(v)}}{n^r}, &\mbox{if $f\in B_{r}(v)$}, \quad r>0,
\end{array}
\right.\qquad\qquad \C\ne\C(n,f).
\end{equation}
\end{lemma}
\subsection{Mapping properties of the Hilbert transform}
The approximation space $B_0(v)$ was firstly introduced in \cite{LutherRusso} in order to find the ``correct  and  minimal'' space for studying the boundedness of the Hilbert transform w.r.t. uniform norm. More precisely, using the standard notation for the constants
\[
c_+:=\max\{0, c\}, \qquad c_-:=\max\{0,-c\},
\]
by \cite[Th. 3.1]{LutherRusso} we have the following
\begin{theorem}\label{th-Hilbert}
	Let $u=v^{a,b}$, with $a,b>-1$, be the Jacobi weight defining the Hilbert transform (\ref{iniziale}) and set $u=\frac{u_+}{u_-}$ with $u_+(x):=v^{a_+,b_+}(x)$ and $u_-(x):=v^{a_-,b_-}(x)$.
	
	Let us consider the following weight functions
	\[
	u_1(x):=u_+(x), \qquad u_2(x):=\frac{u_-(x)}{1+ (u_+u_-)(-1)|\log(1+x)|+ (u_+u_-)(1)|\log(1-x)|}.
	\]

	Hence for all $t\in (-1,1)$ and any $f\in B_0(u_1)$, we have
	\begin{equation}\label{stimaUR}
		|\H^uf(t)|u_2(t)\le \C\left(|f(t)|u_1(t)+\|f\|_{B_0(u_1)}\right), \qquad \C\ne \C(f,t).
	\end{equation}
\end{theorem}
From the previous result it immediately follows that $\H^u$ is a bounded map from $B_0(u_1)$ to $C^0_{u_2}$.

\section{Filtered VP discrete approximation}
This section concerns the main approximation tool we are going to use in getting our quadrature rule for (\ref{iniziale}).

Let $n\in\NN$ and denote by $p_n(w,x)\in \PP_n$ the $n$--th orthonormal polynomial (with  positive leading coefficient) associated to the Jacobi weight
$w(x)=v^{\al,\be}(x)$, with $\al,\be>-1$.
Based on the set of zeros of $p_n(w)$, i. e. $X_n(w)$ in (\ref{Xn}),
the filtered VP polynomial of a function $f$ is defined as follows (see \cite{woula_NA,WCap})
\begin{equation}\label{VP_polynomial}
V_n^m (w,f,x) =\sum_{k=1}^n f(x_k) \Phi_{n,k}^m(x),\qquad 0<m<n,
\end{equation}
with
\begin{equation}\label{fi-filter}
\Phi_{n,k}^m(x)=\lambda_{k}\sum_{j=0}^{n+m-1}\mu_{n,j}^m p_j(w,x)p_j(w,x_k),
\end{equation}
where $\lambda_k$ are the Christoffel numbers defined in (\ref{chris}) and
 $\mu_{n,j}^m$ are the following VP filter coefficients
\begin{equation}\label{muj}
\mu_{n,j}^m:=\left\{\begin{array}{ll}
1, & \mbox{if}\quad j=0,\ldots, n-m,\\ & \\
\displaystyle\frac{n+m-j}{2m}, & \mbox{if}\quad
n-m< j< n+m.
\end{array}\right.
\end{equation}
We point out that the polynomial in (\ref{VP_polynomial}) depends
on two degree--parameters: $n$, which determines the number of nodes, and $m$ which determines the VP filter action range. Both the parameters are involved in the following polynomial preserving property \cite{woula_NA}.
\begin{equation}\label{inva}
V_n^m(w,P)=P,\qquad \forall P\in\PP_{n-m}.
\end{equation}
In what follows we will assume that the degree--parameters $m,n$ are such that
\[
c_1 m\le n\le c_2 m, \qquad \mbox{for some $c_2\geq c_1> 1$ independent of $n$ and $m$},
\]
and we will write $m\approx n$ to express this kind of relation.

We firstly recall the following theorem which gives sufficient conditions such that the map $V_n^m(w):f\in C^0_v\rightarrow V_n^m(w, f)\in C^0_v$  is uniformly bounded w.r.t. $n$ \cite{ThBa}.

\begin{theorem}\label{theo1}
Let $v=v^{\gamma,\delta}$  with $\gamma,\delta\ge 0$ and let $f\in C^0_v$. Moreover let $n\in\NN$ and $m \approx n$ be arbitrarily fixed  and consider  the filtered VP polynomial $V_n^m(w,f)$ associated to the Jacobi weight $w=v^{\alpha,\beta}$, $\alpha,\beta>-1$. If we have
\begin{equation}
\label{theo1ipo1}
\left| \gamma-\delta-\frac{\alpha-\beta}2\right|<1,\qquad
\end{equation}
and
\begin{equation}\label{theo1ipo2}
\left\{
\begin{array}{c}
\displaystyle
\frac{\al}2-\frac 1 4<\gamma\le \frac{\al}2+\frac 5 4 \\ \\
\displaystyle
\frac{\be}2-\frac 1 4<\delta\le \frac{\be}2+\frac 5 4
\end{array}
\right.
\end{equation}
then we get
 \begin{equation} \label{tesi1}
 \|V_n^m (w,f)v\|_\infty\le \C
 \|fv\|_\infty ,\qquad \C\neq \C(n,m,f).
 \end{equation}
\end{theorem}
We also recall that this result has been recently improved in the special case $w$ is one of the four Chebyshev weights, i.e. $|\al|=|\be|=\frac 1 2$, and we refer the reader to \cite[Th.3.1]{Occo_Them} where necessary and sufficient conditions on the weights $v$ and $w$ have been stated.

We point out that, due to (\ref{inva}), the uniform boundedness result (\ref{tesi1}) is equivalent to the following error estimate
\begin{equation}\label{nearbest}
E_{n+m-1}(f)_v\le \|[f-V_n^m(w,f)]v\|_\infty\le \C E_{n-m}(f)_v,\qquad \C\ne \C(n,m,f).
\end{equation}
Hence, taking into account that $m \approx n$ implies $(n-m)\sim n$, we can say that if (\ref{theo1ipo1})--(\ref{theo1ipo2}) hold then we have
\[
\lim_{m \approx n\rightarrow \infty } \|[f-V_n^m(w,f)]v\|_\infty=0, \qquad \forall f\in C^0_v,
\]
being the convergence order comparable with that one of the error of best polynomial approximation $E_{n}(f)_v$.

In the sequel we are going to investigate the behaviour of VP filtered approximation in the case $f$ belongs to the Besov type spaces  $B_r(v)$ introduced in the previous section. We underline that recently in \cite{Mata2020}, some convergence estimates were stated in Zygmund type spaces (i.e. in the Besov type space $B_{r,q}^\infty(v), q=\infty$, $r>0$).

Obviously, under the assumption of Theorem \ref{theo1}, by (\ref{nearbest}) and Lemma \ref{lemmaEm} we easily get
\begin{equation}\label{stimaVP_1}
\|[f-V_n^m (w,f)]v\|_\infty  \leq \C
\left\{ \begin{array}{ll}
\displaystyle \frac{\|f\|_{B_{0}(v)}}{\log n}, & \mbox{if $\ f\in B_{0}(v)$,} \\ &\\
\displaystyle \frac{\|f\|_{B_{r}(v)}}{n^r} ,& \mbox{if $ f\in B_{r}(v)$}, \quad r>0,
\end{array}\right.
\end{equation}	
where $\C\ne \C(n,f)$.

Moreover, we have the following
\begin{theorem} \label{th-VP_Besov}
Let $v=v^{\gamma,\delta}$ and $w=v^{\alpha,\beta}$ be such that  $\gamma,\delta\ge 0$, $\alpha,\beta>-1$ and assume $V_n^m(w):C^0_{v}\rightarrow C^0_{v}$ is an uniformly bounded map w.r.t $n,m\in\NN$ with $m \approx n$. For any  $r\geq0$, also the map $V_n^m(w): B_{r}(v)\rightarrow B_{r}(v)$ is uniformly bounded w.r.t. $n$ and we have
\begin{equation}\label{eq-conv}
\lim_{m \approx n\rightarrow\infty}\|f-V_n^m(w,f)\|_{B_r(v)}=0, \qquad \forall f\in B_r(v), \quad \forall r\ge 0.
\end{equation}
Moreover, for all $f\in B_{r}(v)$ with $r>0$, the following error estimates hold
\begin{eqnarray}\label{VP_map1}
\|f-V_n^m  (w,f)\|_{B_{0}(v)}&\le & \C\frac{\log n}{n^{r}}\ \|f\|_{B_{r}(v)},\qquad
r>0,\\
\|f-V_n^m  (w,f)\|_{B_{s}(v)}&\le &\frac{\C}{n^{r-s}}\ \|f\|_{B_{r}(v)}, \qquad r\ge s>0,\label{VP_map2}
\end{eqnarray}
where in both cases $\C\ne \C(n,f)$.
\end{theorem}
\subsection{The limit case $m=0$:  the Lagrange interpolation}

As said in the Introduction, for all $n\in \NN$, when the parameter $m$ defining the VP operator is set equal to zero, we get  the Lagrange interpolating operator at the same system of nodes $X_n(w)$,  namely
\begin{equation*}
	L_n(w,f,x)=\sum_{k=1}^n f(x_k) l_{n,k}(x), \qquad l_{n,k}(x)=\lambda_k\sum_{j=0}^{n-1}p_j(w,x_k)p_j(w,x),
\end{equation*}	
being $\lambda_k$ given in (\ref{chris}).

While under suitable conditions for the weights $w$ and $v$  the VP operator $V_n^m(w):C^0_v\rightarrow C^0_v$ is uniformly bounded w.r.t. $n$  (cf. Th. \ref{theo1}), it is well known that the norm of the Lagrange operator $L_n(w):C^0_v\rightarrow C^0_v$, the so called weighted Lebesgue constant, does not. Nevertheless, necessary and sufficient conditions for optimal Lebesgue constant, i. e.  behaving like $\log n$, are known \cite{MR_Besov_inf,mastromilobook}.
More precisely, we have the following theorem which collects the analogous of previous estimates stated for the VP approximation.
\begin{theorem}\label{convLag}
Let $v=v^{\gamma,\delta}$ and $w=v^{\alpha,\beta}$ be such that  $\gamma,\delta\ge 0$, $\alpha,\beta>-1$. For all $f\in C^0_v$, the conditions
\begin{equation}\label{condConvLag}
	\left\{
	\begin{array}{c}
		\displaystyle
		\frac{\al}2+\frac 1 4\le\gamma\le \frac{\al}2+\frac 5 4 \\ \\
		\displaystyle
		\frac{\be}2+\frac 1 4\le\delta\le \frac{\be}2+\frac 5 4
	\end{array}
	\right.
\end{equation}
are necessary and sufficient for having
\begin{equation}\label{lagrange}
	\|[f-L_n(w,f)]v\|_\infty \leq \C \log n E_n(f)_v, \quad \C\neq \C(n,f).
\end{equation}
Moreover, if (\ref{condConvLag}) holds then for all $f\in B_{r}(v)$ with $r>0$, we get
\begin{eqnarray}\label{Lag_map1}
	\|f-L_n(w,f)\|_{B_{0}(v)}&\le & \C\frac{\log^2 n}{n^{r}}\ \|f\|_{B_{r}(v)},\qquad
	r>0,\\
	\|f-L_n(w,f)\|_{B_{s}(v)}&\le &\C \frac{\log n}{n^{r-s}}\ \|f\|_{B_{r}(v)}, \qquad r\ge s>0,\label{Lag_map2}
\end{eqnarray}
where in both cases $\C\ne \C(n,f)$.
\begin{remark}
Comparing the convergence results between VP and Lagrange approximations we remark that for the same class of function (say $f\in B_r(v)$, $r>0$) the convergence rate is similar, but in the Lagrange estimates an extra factor $\log n$ appears. Moreover conditions (\ref{theo1ipo2}) are wider than (\ref{condConvLag}).
\end{remark}	
\end{theorem}
\section{Filtered VP quadrature rules for the Hilbert transform}
In this section we are going to introduce a new class of product integration rules for the Hilbert transform in (\ref{iniziale}) based on filtered VP approximation. More precisely,  by
replacing  $f$ in (\ref{iniziale})  with the polynomial $ V_n^m (w,f)$, we define
\begin{equation}\label{def-gen}
	\mathcal{H}^u_{n,m}(w,f,t):=\int_{-1}^{1} \frac {V_n^m (w,f,x)}{x-t}u(x)dx, \qquad -1<t<1.
\end{equation}
More explicitly, by (\ref{VP_polynomial})-(\ref{fi-filter}) we get the following quadrature rule (VP rule)
\begin{equation*}\label{def-compact}
	\mathcal{H}^u_{n,m}(w,f,t)=\sum_{j=0}^{n+m-1}\rho_{n,j}^m(w,f) Q_j^u(w,t),\qquad -1<t<1,
\end{equation*}
where we set
\begin{eqnarray}
	\label{Q}
	Q_j^u(w,t)&:=&\int_{-1}^{1}\frac {p_j(w,x)}{x-t}u(x)dx,\\
	\rho_{n,j}^m(w,f)&:=&\mu_{n,j}^m \sum_{k=1}^n\lambda_k p_j(w,x_k)f(x_k), \nonumber
\end{eqnarray}
and with $\mu_{n,j}^m$ defined in (\ref{muj}).
\subsection{Computational details}
First of all we point out that the coefficients $\rho^m_{n,j}(w,f)$ do not depend on $t$, and therefore are not influenced by the closeness of $t$ to any interpolation point $\{x_i\}_{i=1}^n$.
This means that the numerical cancellation is avoided in computing these quantities.

About  the functions  $\{Q_j^u(w,t)\}_j$, we remark that only in some particular cases they are explicitly known  (see e.g. \cite{PS}) while in the general case
they can be computed via a recurrence relation (see \cite{king}). The following proposition states such relation, deducible from the well--known recurrence relation for the orthonormal sequence $\{p_j(w)\}_{j}$
\begin{equation}\label{3term-pol}
	b_{j+1} p_{j+1}(w,x)=(x-a_{j})p_{j}(w,x)- b_j p_{j-1}(w,x),\qquad j\ge 0,
\end{equation}
being $p_{-1}(w,x)=0$ and $p_0(w,x)=\left(\int_{-1}^1 w(x)dx\right)^{-\frac 12}$ the starting values (see for instance \cite[pp.131-133]{mastromilobook} for the explicit expressions of the coefficients).
\begin{proposition}\label{th-recQ}
	For all Jacobi weights $u,w$ and any $t\in (-1,1)$, the functions $Q_j(t):=Q_j^u(w,t)$ defined in (\ref{Q}) satisfy the following three-term recurrence relation
	\[
	Q_{j+1}(t)=(A_jt+B_{j})Q_{j}(t)- C_j Q_{j-1}(t)+D_j,\qquad  j\ge 0,
	\]
	where the starting values are given by
	\begin{equation*}\label{recurrQ}
		Q_{-1}(t)=0, \qquad Q_0(t)=\left(\int_{-1}^1\frac{u(x)}{x-t}dt\right)\left(\int_{-1}^1 w(t)dt\right)^{-\frac 12},
	\end{equation*}
	and the coefficients are defined, by means of the coefficients $\{a_j\}_j$ and $\{b_j\}_j$ in (\ref{3term-pol}), as follows
	\[
	A_j=b_{j+1}^{-1},\ B_j=-\frac{a_j}{b_{j+1}},\ C_j=\frac{b_j}{b_{j+1}},
	\ D_j=\frac{1}{b_{j+1}} \int_{-1}^1 p_j(w,x) u(x)dx,\quad j=0,1,\dots.
	\]
\end{proposition}

\subsection{Error estimates}
Let us analyze the error function
\begin{equation}\label{def-err}
	\E_{n,m}^u(w,f,t):=\left|\H^uf(t)-\mathcal{H}_{n,m}^u(w,f,t)\right|,
	\qquad -1<t<1.
\end{equation}
Of course the behaviour of such error is influenced by the approximation error provided by the VP polynomial that has been employed. Such dependence is specified in the following convergence theorem
\begin{theorem}\label{theo2}
	Under the notation of Theorem \ref{th-Hilbert}, for all $f\in B_0(u_1)$ and any $t\in(-1,1)$, we have
\begin{equation}\label{errore_quad_point}
		\E_{n,m}^u(w,f,t)u_2(t) \le \C \left|f(t)-V_n^m(w,f,t)\right| u_1(t)	
+  \C \|f-V_n^m (w, f)\|_{B_{0}(u_1)},
	\end{equation}
	where $\C\ne \C(n,f,t)$.
Moreover, if $V_n^m(w):C^0_{u_1}\rightarrow C^0_{u_1}$ is an uniformly bounded map w.r.t. $m \approx n$, then
	\begin{equation}\label{conv-B0}
		\lim_{n\rightarrow \infty}\E_{n,m}^u(w,f,t)u_2(t)=0, \qquad \forall f\in B_0(u_1),
	\end{equation}
	and the convergence holds uniformly in $t\in (-1,1)$.
\end{theorem}
The previous theorem assures that if the function $f$ is in $B_{0}(u_1)$ then the proposed quadrature formula converges on conditions that the weight $w=v^{\alpha,\beta}$ is suitably chosen.
In particular Theorem \ref{theo2} and Theorem \ref{theo1} provide the following criteria for the choice of the weight $w=v^{\alpha,\beta}$
\begin{equation}\label{criterio}
	 a_+-b_+-1<\frac{\alpha-\beta}2<a_+-b_++1,\qquad
	\mbox{and}\qquad
	\left\{
	\begin{array}{c}
		\displaystyle
	2 a_+-\frac 5 2\leq \alpha < 2 a_+ +\frac 1 2 \\ \\
		\displaystyle
			2 b_+-\frac 5 2\leq \beta < 2 b_+ +\frac 1 2
	\end{array}
	\right.
\end{equation}
where $a_+=\max\{a,0\}$,  $b_+=\max\{b,0\}$, being $u=v^{a,b}$ the weight defining the Hilbert transform.

Moreover, recalling Theorem \ref{th-VP_Besov}, under the hypotheses of Theorem \ref{theo2} the map $V_n^m(w):B_0(u_1)\rightarrow B_0(u_1)$ is uniformly bounded w.r.t. $n\sim (n-m)$. Hence by the invariance property (\ref{inva}) and by (\ref{errore_quad_point}) we get
\begin{equation}\label{err-ptuale}
	\E_{n,m}^u(w,f,t)u_2(t) \le \C |f(t)-V_n^m(w,f,t)|u_1(t)	+ \C E_n(f)_{B_{0}(u_1)},
\end{equation}
where $  E_n(f)_{B_{0}(u_1)}:=\inf_{P\in\PP_n}\|f-P\|_{B_{0}(u_1)}$ and $\C\neq\C(n,m,f)$.

We point out that (\ref{err-ptuale}) says that for each $t$ the convergence rate of the VP-rule depends on two components: the pointwise approximation that the VP polynomial of $f$ provides at the specific $t\in (-1,1)$, and the smoothness of $f$ that influences the convergence rate of the error of best polynomial approximation of $f$.

If $f$ is smoother than in Theorem \ref{theo2}, from  (\ref{errore_quad_point}) and (\ref{VP_map1}) we immediately deduce the following corollary.
\begin{corollary}\label{corollarioConv}
	Under the same assumptions of Theorem \ref{theo2}, for all $f\in B_r(u_1)$, $r>0$, we get
	\begin{equation}
		\label{errore_quad2}
		\E_{n,m}^u(w,f,t)u_2(t) \le \C  |f(t)-V_n^m(w,f,t)|u_1(t)	+\C \frac{\log n}{n^r} \|f\|_{B_{r}(u_1)},
	\end{equation}
where $\C\neq \C(n,f,t).$
\end{corollary}
Obviously, concerning the converge rate of the VP-rule, from (\ref{errore_quad2}) and (\ref{stimaVP_1}) we get that
\begin{equation}
	\label{errore_quad}
	\E_{n,m}^u(w,f,t)u_2(t) \le \C\frac{\log n}{n^{r}}, \qquad \C\ne \C(n,t),
\end{equation}
holds true $\forall f\in B_r(u_1)$, with $r>0$.

We conclude the section with a comparison between the proposed VP-rule and the analogous one using the Lagrange polynomial (L-rule):
\begin{equation}\label{prodLag}
	\H^u_nf(t):=\int_{-1}^1 \frac{L_n(w,f,x)}{x-t}u(t)\ dx.
\end{equation}
This kind of formula was considered by several authors in the case $w=u$ (see for instance \cite{CMR} and the reference therein). Here, set
\begin{equation*}\label{def-err-Lag}
	\E_{n}^u(w,f,t):=\left|\H^uf(t)-\mathcal{H}_{n}^u(w,f,t)\right|,
	\qquad -1<t<1,
\end{equation*}
as a consequence of Theorems \ref{th-Hilbert} and \ref{convLag}, we note that if the exponents of $u=v^{a,b}$ and $ w=v^{\alpha,\beta}$ are such that
\begin{equation}\label{criterio-L}
\left \{\begin{array}{c}
	\displaystyle
2 a_+-\frac 5 2\leq \alpha \leq 2 a_+ -\frac 1 2 \\ \\
\displaystyle
2 b_+-\frac 5 2\leq \beta\leq 2 b_+ -\frac 1 2\end{array} \right.
\end{equation}
then, for all $f\in B_r(u_1)$, $r>0$, we have
	\begin{equation}
		\label{errore_quad3}
		\E_{n}^u(w,f,t)u_2(t) \le \C |f(t)-L_n(w,f,t)|u_1(t)	+\C\frac{\log^2 n}{n^{r}}\ \|f\|_{B_{r}(u_1)},
	\end{equation}
with $\C\ne \C(n,f,t).$
\section{Numerical experiments}
In this section  we report  some numerical tests in order to evaluate the performance of the proposed VP--rule, also in comparison with other quadrature formulae.

To be more precise, for several choices of $t\in (-1,1)$, we  compare the absolute errors $e^{VP}$ achieved  by our  VP-rules,  with the errors $e^L$ obtained by the L-rule in (\ref{prodLag}).  Moreover, we will make a comparison with the Modified Gaussian rule  (MG-rule) proposed in \cite{crimas_gau,CMR}. Finally, we will also consider an example given in \cite{hasegawa} by Hasegawa and Torii, who proposed a product rule (HT--rule) whose absolute error we denote by $e^{HT}$.

Since the exact values are not available in the general case, we have retained ``exact''  the values achieved  by means of higher degree quadrature formulas, varying the choice of the rule among the previous  VP, L and MG rules. In each test we will declare the rule we have used to determine  the exact value.

From the wide experimentation we carried on, we have selected five examples, varying the regularity of the  functions $f$ and the possible choices of the weight functions  $u$ and $w$, which define the transform $\mathcal{H}^u(f)$ and the VP  approximants $V_n^m(w,f)$, respectively.

In all the tests   $n$ will denote the number of nodes employed in the considered rules.
Only  the MG-rule  makes exception, since it is based on a variable number of nodes that could be $n$ or $n+1$, being this choice contemplated by the method itself, according to the position of $t$ \cite{CMR}.

The numerical outputs of each example are collected in tables having the same numbering as the examples to which they refer. At each row of the tables, besides the number of nodes $n$, the value of the additional parameter $m$  chosen in the VP-rule  is specified and the best quadrature error is evidenced in bold.

Finally,
at the end of the section we present two further tests on the VP-rules, where we focus on the behaviour of the absolute errors  as $n$ increases, making several choices of $m\approx n$. More precisely, in Test A we take $m=\lfloor \theta n\rfloor$ with fixed $\theta\in (0,1)$ and we compare the error behaviours corresponding to different choices of $\theta$. In Test B, in correspondence of any fixed number of nodes we plotted the behaviour of the best and the worst $e^{VP}$ obtained for $0<m<n$.

In the sequel details and comments concerning the five examples and the two tests are given. We point out that all the computations have been performed in double-machine precision ($eps_D \approx 2.22044e-16$) by using MATLAB R2021a.

\noindent{\bf Example 1}
\[\mathcal{H}^uf(t)=\int_{-1}^1\frac{|x-0.5|^{10.01}}{x-t }\sqrt{\frac{1-x}{1+x}}dx. \]
Here $ u=v^{\frac 1 2 ,-\frac 1 2 }$, $u_1=v^{\frac 1 2, 0}$, $u_2=v^{0,\frac 1 2}$ and we chose $w=u$. This choice ensures that (\ref{criterio}) and (\ref{criterio-L}) are satisfied (with $\alpha=-\beta=\frac 1 2$, $a_+=\frac 12$ and $b_+=0$).
The function $f\in B_r(u_1),$ with $r<10.01$ and according to (\ref{errore_quad}), the  error goes like $\mathcal{O}\left(\frac{\log n}{n^{10.01}}\right)$. Taking into account that the seminorm
$\|f\|_{B_r(u_1)}\ge 2.5\times 10^6$,   the theoretical rate of convergence  is confirmed, also for values  of $t$ ``close''
to the point $0.5$ where the function $f$ is less regular. Indeed, for $n=151$ we can expect $14$ exact digits.

In order to compute the absolute errors,  we have retained exact the values achieved  by means of the MG rule of order $n=1200.$
Moreover since in this case $w=u=v^{a,-a}$, with $0<a<1$, for the computation of the coefficient of the product rules we have used the following
explicit expression (see \cite[p.310]{PS})
\begin{equation*}\label{eq-Qspecial}
	Q_j^u(u,t)=\int_{-1}^1\frac{p_j(u,x)}{x-t }u(x)dx=\pi\cot (\pi a) p_j(u,t)u(t)-\frac\pi{\sin(\pi a)}p_j(u^{-1},t),
	\quad j\ge 0.
\end{equation*}

\begin{table}[h]
	\centering
{\footnotesize 	\begin{tabular}{|c|c|l|l|l|c|c|l|l|l|}
		\hline
		\multicolumn{5}{|c|}{$\mathbf{t=0.499999999}$} & \multicolumn{5}{c|}{$\mathbf{t=0.5}$}\\
		\hline\hline
		$n$    &  $m$    &  $e_n^{VP},$   &  $e_n^{L}$   &   $e_n^{MG}$  & $n$    &  $m$    &  $e_n^{VP},$   &  $e_n^{L}$   &   $e_n^{MG}$ \\ \hline
		51    &8	     &  8.53e-14	 & 6.47e-13	 & \textbf{2.84e-14}  & 51	 &8	     &  9.24e-14	 & 6.61e-13	 & \textbf{2.13e-14}\\ \hline
		151	 &13	 &  \textbf{0.00}	 & 1.22e-12	 & 1.42e-14  & 151	 &19	 & \textbf{0.00}	 & 1.23e-12	 & 2.12e-14\\ \hline\hline
		\multicolumn{5}{|c|}{$\mathbf{t=0.51111}$} & \multicolumn{5}{c|}{$\mathbf{t=0.75}$}\\
		\hline\hline
		$n$    &  $m$    &  $e_n^{VP},$   &  $e_n^{L}$   &   $e_n^{MG}$  & $n$    &  $m$    &  $e_n^{VP},$   &  $e_n^{L}$   &   $e_n^{MG}$ \\ \hline
		51	 &4	 &  6.39e-14	 & 5.12e-13	 & \textbf{2.84e-14} &  51	     &3	     &  \textbf{0.00}	 & 4.33e-13	 & 4.62e-14 \\ \hline
		151	 &68	 &  \textbf{0.00}	 & 1.61e-12	 & 8.53e-14 & 	 &	 & 	 & 	 & \\ \hline
	\end{tabular}}
	\caption{Example 1 }
\end{table}

\noindent{\bf Example 2}
\[\mathcal{H}^u f(t)=\int_{-1}^1\frac{{\log (1-x)}}{x-t }(1-x)^{0.4}(1+x)^{0.25}dx. \]
Here $u\equiv u_1=v^{0.4 ,0.25}$, $u_2\equiv 1$ and we chose $w=v^{-\frac 1 2 ,-\frac 1 2 }$, whose exponents satisfy (\ref{criterio}) and (\ref{criterio-L}) with $a_+=0.4$ and $b_+=0.25$

In this case
$\ f\in B_{r}(u),\ r<0.8$, and therefore, by (\ref{errore_quad}) the error of the VP-rule is bounded by  $\mathcal{O}\left(\frac{\log n}{n^{r}}\right)$.
Also in this test,  we have retained exact the values achieved  by means of the MG rule of order $n=1200$ and the results are shown in Table 2.

For $t$ approching to $1$,  we observe a progressive  loss of exact digits, coherent with the error estimate, while the situation appears much more better when $t$ lies in the remaining part of the interval. Overall, the VP-rule seems to provide a good performance, on average better than those offered by the other two rules.

\begin{table}[h]
	\centering
{\footnotesize 	\begin{tabular}{|c|c|l|l|l|c|c|l|l|l|}
		\hline
		\multicolumn{5}{|c|}{$\mathbf{t=-0.1}$} & \multicolumn{5}{c|}{$\mathbf{t=0.7}$}\\
		\hline
		$n$    &  $m$    &  $e_n^{VP},$   &  $e_n^{L}$   &   $e_n^{MG}$  & $n$    &  $m$    &  $e_n^{VP},$   &  $e_n^{L}$   &   $e_n^{MG}$ \\ \hline
		50   & 18 &  2.59e-05  &   9.57e-05   &  \textbf{6.01e-06 } & 50	 &30	 &  \textbf{4.11e-07}	 & 2.71e-04	 & 2.20e-05\\ \hline
		101  & 78 &   4.25e-06  &   3.78e-05  &   \textbf{1.45e-06 } & 150	 &120	 &  \textbf{3.98e-11}	 & 6.77e-06	 & 2.12e-06\\ \hline
		201  & 78 &    \textbf{4.81e-08 } &   4.24e-06  &   2.92e-07 & 250	 &14	 &  7.74e-06	 & 2.65e-05	 & \textbf{6.29e-07}\\ \hline
		301  &31&   \textbf{7.91e-08 } &   2.17e-06   &  1.08e-07 & 350	 &156	 &  4.30e-07	 & 1.62e-05	 & \textbf{2.67e-07} \\ \hline
		400 &377&    \textbf{1.61e-08}  &   3.00e-06   &  5.18e-08 & 450	 &164	 & \textbf{1.24e-07}	 & 2.37e-06	 & 1.38e-07\\ \hline
		500  &313&   \textbf{1.53e-08}    & 5.32e-07   &  2.83e-08 & 550	 &389	 & \textbf{7.91e-08}	 & 5.24e-06	 & 7.98e-08\\ \hline\hline
	\end{tabular}}
\end{table}

\begin{table}[h]
	\centering
{\footnotesize 	\begin{tabular}{|c|c|l|l|l|c|c|l|l|l|}
		\hline
		\multicolumn{5}{|c|}{$\mathbf{t=0.9}$} & \multicolumn{5}{c|}{$\mathbf{t=0.9999}$}\\
		\hline
		$n$    &  $m$    &  $e_n^{VP},$   &  $e_n^{L}$   &   $e_n^{MG}$  & $n$    &  $m$    &  $e_n^{VP},$   &  $e_n^{L}$   &   $e_n^{MG}$ \\ \hline
		50	 &10	 &  \textbf{5.74e-05}	 & 4.22e-04	 & 6.60e-05 & 50	 &3	 &  \textbf{1.78e-02}	 & 1.78e-02	 & 3.40e-02 \\ \hline
		150	 &45	 & \textbf{8.24e-07}	 & 1.90e-04	 & 6.45e-06 & 250	 &236	 & {4.25e-03}	 & 5.08e-03 & \textbf{1.65e-03} \\ \hline
		250	 &150	 &  4.50e-06	 & 2.66e-05	 & \textbf{1.87e-06} & 350	 &271	 &  \textbf{4.96e-05}	 & 6.53e-03	 & 7.51e-04\\ \hline
		350	 &38	 &  \textbf{1.66e-07}	 & 2.82e-05	 & 8.08e-07 & 550	 &495	 &  \textbf{5.40e-05}	 & 2.63e-03	 & 2.28e-04\\ \hline
		450	 &66	 &  \textbf{8.47e-08}	 & 2.36e-05	 & 4.16e-07 &  650	 &525	 &  \textbf{5.50e-06}	 & 3.40e-04	 & 1.37e-04\\ \hline \hline
	\end{tabular}}
	\caption{Example 2}
\end{table}

\noindent{\bf Example 3}
\[\mathcal{H}(f,t)=\int_{-1}^1\frac{f(x)}{x-t }dx, \]
$$f(x)= e^{8(x-1)},\ u=v^{0,0}=u_1,  \ f\in B_r, \ \forall r>0,\  w=v^{-\frac 1 2,-\frac 1 2}.$$
This example is taken from \cite{hasegawa} where the quadrature rule here denoted by HT--rule has been proposed. Hence, in Table 3 we have an additional column reporting the errors $e^{HT}$ of such rule as provided in \cite{hasegawa}, writing n.a. (not available) in the cases they are not furnished by the authors.
\begin{table}[h]
	\centering
{\footnotesize  	\begin{tabular}{|c|c|l|l|l|l|c|c|l|l|l|l|}
		\hline
		\multicolumn{6}{|c|}{$\mathbf{t=0.2}$} & \multicolumn{6}{|c|}{$\mathbf{t=0.5}$}\\
		\hline\hline
		$n$    &  $m$    &  $e_n^{VP}$   &  $e_n^{L}$   &   $e_n^{MG}$  & $e_n^{HT}$ & $n$    &  $m$    &  $e_n^{VP},$   &  $e_n^{L}$   &   $e_n^{MG}$  & $e_n^{HT}$\\ \hline
		25&     2&     1.4e-13 &    8.0e-14  &   \textbf{1.4e-15} &  9.0e-13 &  25	 &2	 &  6.5e-13	  	 & 2.9e-13  & \textbf{2.5e-15}& 3.0e-13\\ \hline
		30	 &2	 & \textbf{1.9e-15}	 & 4.2e-15	 & 2.7e-15 & n.a.&  51	 & 27	 & \textbf{3.8e-16}	 & 3.6e-14	 & 5.8e-15 & n.a.\\ \hline
		101	 &60	 &  \textbf{5.6e-17}	 & {8.1e-14}	 & 2.9e-15 & n.a. &&&&& &\\ \hline \hline
		\multicolumn{6}{|c|}{$\mathbf{t=0.95}$} & \multicolumn{6}{|c|}{$\mathbf{t=0.999}$}\\
		\hline\hline
		$n$ &$m$&  $e_n^{VP}$       & $e_n^{L}$   &$e_n^{MG}$        & $e_n^{HT}$ & $n$    &  $m$   &  $e_n^{VP}$   &  $e_n^{L}$   &   $e_n^{MG}$  & $e_n^{HT}$ \\ \hline
		25	 &2	 &  1.2e-12	         & 3.5e-13	   &\textbf{1.3e-14}  & 1.0e-12    & 25 	&   2	 &  8.0e-13	                             & 1.9e-13  &  \textbf{1.8e-15} & n.a.\\ \hline
		30	 &2	 & \textbf{3.4e-15}	 & 1.6e-14     & 1.4e-14          &      n.a.      & 30     &2	     & \textbf{8.9e-16}		 & 1.8e-15 & {1.1e-14} & n.a.\\ \hline
		51	 &3	 &   \textbf{3.1e-15}& 2.2e-14     & 2.5e-14 & n.a. &&&&&&\\ \hline\hline
	\end{tabular}}
	\caption{Example 3}
\end{table}

The function $f$ is very smooth and, as we expect,  the VP, L and  MG rule give a good performance. We underline that speeding up a little bit $n$ the VP-rule catches the machine precision.

\noindent{\bf Example 4}
\[\mathcal{H}^u(f,t)=\int_{-1}^1 \frac 1 {x^2+2^{-10}} \frac{\sqrt[3]{1-x^2}}{x-t }dx, \]

Here $ u\equiv u_1=v^{\frac 1 3 ,\frac 1 3}$, $u_2\equiv 1$ while we fixed $w=v^{1 ,1}.$

As already remarked the sufficient conditions stated for the convergence of the VP-rule are wider than those of the L-rule. In this case $a_+=b_+=\frac 1 3$ and the choice of $w$ assures the convergence of the VP-rule, but not that of the L-rule.
So we will compare the results by VP-rules only with those by MG-rule.
In this test,  we have retained exact the values achieved  by means of the VP-rule of order $n=1000, m=500.$

In this example the function $\ f\in B_r(u), \ \forall r>0,$ but  the presence of the complex conjugate  poles $\pm i 2^{-5}$ too close to the integration interval $(-1,1)$ produces slower convergence, as well as happens in quadrature rules for ordinary integrals (see e.g. \cite{davis,Monegato_poles}).
The numerical results in Table 4 confirm this trend for the
MG--rule that shows a certain saturation, while it seems that the VP-rule converges faster.

\begin{table}[h]
	\centering
{\footnotesize 	\begin{tabular}{|c|c|l|l|c|c|l|l|}
		\hline
		\multicolumn{4}{|c|}{$\mathbf{t=0.2}$} & \multicolumn{4}{c|}{$\mathbf{t=0.4}$}\\
		\hline
		$n$    &  $m$    &  $e_n^{VP},$    &   $e_n^{MG}$  & $n$    &  $m$    &  $e_n^{VP}$     &   $e_n^{MG}$ \\ \hline
		81	 &48	 & \textbf{7.14e-02}	 & 5.94e+00	 & 81	 &40	 &  \textbf{8.69e-02}	 & 2.81e+00	\\ \hline
		101	 &90	 & \textbf{ 1.58e-02}	 & 1.69e+00	 & 101	 &9	 &  \textbf{8.87e-03}	 & 8.63e-01\\ \hline
		201	 &61	 &  7.08e-03	 & \textbf{3.06e-03}	 &  201	 &160	 &  \textbf{5.93e-04}	 & 1.66e-03\\ \hline
		301	 &30	 &  6.69e-05	 & \textbf{5.92e-06}	 &  401	 &160	 &  3.31e-05	 & \textbf{2.10e-08}\\ \hline
		401	 &30	 &  7.11e-05	 & \textbf{4.08e-07}	 &  501	 &53	 &  4.06e-07	 & \textbf{1.47e-08}\\ \hline
		501	 &30	 &  \textbf{7.06e-08}	 & 3.96e-07	 &  601	 &38	 &  \textbf{7.67e-09}	 & {1.48e-08}	\\ \hline
		601	 &34	 &  \textbf{6.36e-10}	 & 3.96e-07	 &  701	 &51	 &  \textbf{3.57e-10}	 & 1.47e-08	\\ \hline \hline
		\multicolumn{4}{|c|}{$\mathbf{t=0.6}$} & \multicolumn{4}{c|}{$\mathbf{t=0.9}$}\\
		\hline\hline
		$n$    &  $m$    &  $e_n^{VP},$    &   $e_n^{MG}$  & $n$    &  $m$    &  $e_n^{VP}$     &   $e_n^{MG}$ \\ \hline
		81	 &55	 & \textbf{6.98e-02}	 & 1.88e+00	 & 81	 &71	 &  \textbf{7.20e-01}	 & 1.35e+00	 \\ \hline
		101	 &9	 &  8.71e-01	 &\textbf{ 5.40e-01}	 &  101	 &9	 &  \textbf{2.09e-01}	 & 3.85e-01	 \\ \hline
		201	 &19	 &  \textbf{3.78e-04}	 & 1.11e-03	 &    201	 &28	 &  \textbf{6.28e-04}	 & 6.99e-04	                          \\ \hline
		301	 &34	 &  5.36e-05	 & \textbf{2.01e-06	} &    301	 &84	 &  1.10e-04	 & \textbf{1.36e-06 }     \\ \hline
		401	 &200	 & {8.17e-07}	 & \textbf{1.21e-08}	 &    401	 &28	 &  {1.30-08}	 & \textbf{3.06e-09}  \\ \hline
		501	 &200	 & {1.84e-07}	 & \textbf{7.98e-09}	 &    501	 &28	 &  \textbf{9.67e-10}	 & 5.84e-09	  \\ \hline
		601	 &83	 &  \textbf{4.91e-10}	 & 7.98e-09	 &   601	 &70	 &  \textbf{1.64e-09}	 & 5.84e-09	    \\ \hline
		701	 &79	 &  \textbf{5.23e-10}	 & 7.98e-09	 &    701	 &144	 &  \textbf{8.73e-12}	 & 5.84e-09	   \\ \hline
		\hline
	\end{tabular}}
	\caption{Example 4 }
\end{table}

\noindent{\bf Example 5}
\[\mathcal{H}^uf(t)=\int_{-1}^1\frac{f(x)}{x-t }\sqrt{1-x^2}dx, \]
$$f(x)=\frac 1{1+1000(x+0.5)^2}+\frac 1{\sqrt{1+1000(x-0.5)^2}}.$$
In this example  $ u\equiv u_1=v^{\frac 1 2 ,\frac 1 2}$, $u_2\equiv 1$, and we fixed $ w=v^{-\frac 12 ,-\frac 12}.$ The density function $f$ belongs to $B_r(u_1)$, for any $r>0.$
Nevertheless the graphic of $f$, given in Figure 1, shows two picks corresponding to $x=\pm 0.5$. In this case the pointwise approximation provided by the VP polynomial is almost everywhere better than the Lagrange approximation and we aim to test if such difference is reflected in the quadrature errors $e^{VP}$ and $e^L$ too. For this reason we exclude MG--rule from this test, limiting the comparison between VP and L rules.
\begin{figure}[h]
	\begin{center}
		\includegraphics[scale=.30]{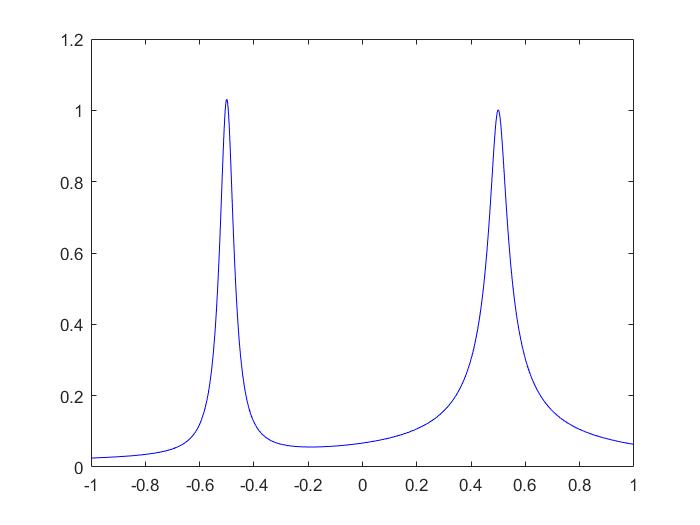}
		{\caption{\label{funzione}Example 5: graphic of $f$ }}
	\end{center}
\end{figure}

Retaining exact the values achieved  by means of the L-rule of order $n=1200$, the errors $e^{VP}$ and $e^L$ are shown in Table 5 for different choices of  $t>0$.  We can observe that  for $t=0.5$ the errors by VP-rules are comparable with those by the L-rule, while in points $t$ far from the pathological points $\pm 0.5$, VP-rule gives better results.
\begin{table}[h]
	\centering
{\footnotesize 	\begin{tabular}{|c|c|l|l|c|c|l|l|}
		\hline
		\multicolumn{4}{|c|}{$\mathbf{t=0.1}$} & \multicolumn{4}{c|}{$\mathbf{t=0.2}$}\\
		\hline
		$n$    &  $m$    &  $e_n^{VP},$ & $e_n^{L}$   & $n$    &  $m$    &  $e_n^{VP}$  & $e_n^{L}$   \\ \hline
		20	 & 12	 &  \textbf{1.79e-03}	 & 2.38e-01	    &       20	 &8	 &  \textbf{5.67e-03}	 & 2.18e-01	     \\ \hline
		30	 &3	     &  \textbf{2.82e-03}	 & 3.37e-03     &      30	 &9	 &  \textbf{2.97e-04}	 & 1.69e-03	     \\ \hline
		40	 &17	 &  \textbf{5.90e-04}	 & 9.98e-02	    &     40	 &19	 &  \textbf{3.90e-03}	 & 1.40e-01	  \\ \hline
		50	 &16	 &  \textbf{5.63e-04}	 & 7.40e-02 	&    50	 &8	 &  \textbf{5.15e-04}	 & 4.84e-02	       \\ \hline
		60	 &6	     & \textbf{1.73e-05}	 & 1.53e-04     &   60	 &24	 &  \textbf{7.65e-05}	 & 3.28e-03	      \\ \hline
		70	 &14	 &  \textbf{3.94e-05}	 & 2.62e-02	    &     70	 &49	 &  \textbf{2.96e-04}	 & 4.17e-02	 \\  \hline
		200 &45	 &  \textbf{1.73e-06}	 & 2.32e-04	         & 200	 &35	 &  \textbf{2.31e-07}	 & 1.32e-04	    \\  \hline
		250 &175	 & \textbf{1.24e-09}	 & 4.23e-06	                 &250	 &24	 &  \textbf{7.45e-08}	 & 3.37e-06	  \\ \hline
		300 &17	 &  \textbf{1.33e-08}	 & 1.36e-06	           &300	 &181	 &  \textbf{3.74e-07}	 & 9.28e-07	    \\  \hline
		\hline
		\multicolumn{4}{|c|}{$\mathbf{t=0.5}$} & \multicolumn{4}{c|}{$\mathbf{t=0.8}$}\\
		\hline
		$n$    &  $m$    &  $e_n^{VP},$ & $e_n^{L}$               & $n$    &  $m$    &  $e_n^{VP}$  & $e_n^{L}$   \\ \hline
		30	 &5	 &  \textbf{3.10e-02}	 & \textbf{3.10e-02}	    &   20	 &2	 &  \textbf{1.11e-02}	 & 1.27e-02	            \\ \hline
		50	 &3	 &  1.42e-01	 & \textbf{1.25e-01}	            &     30	 &24	 &  \textbf{3.75e-02}	 & 3.86e-02	                   \\ \hline
		60	 &23	 &  \textbf{3.92e-03}	 & \textbf{3.92e-03}	&      50	 &20	 &  \textbf{1.82e-04}	 & 9.12e-03	             \\ \hline
		80	 &7	 &  4.86e-02	 & \textbf{4.31e-02}	            &     70	 &13	 &  \textbf{5.37e-05}	 & 6.01e-03	                \\ \hline
		90	 &80	 &  \textbf{4.61e-04}	 & \textbf{4.61e-04}	&    	 80	 &8	 &  \textbf{5.78e-05}	 & 3.79e-03	    \\  \hline
		150 &105	 &  \textbf{6.20e-06}	 & \textbf{6.20e-06}    & 100	 &60	 &  \textbf{1.40e-04}	 & 1.42e-03	                                                           \\ \hline
		200 &13	 &  5.06e-04	 & \textbf{4.19e-04}	      &         150	 &11	 &  \textbf{3.41e-06}	 & 3.82e-05	                                                   \\ \hline
		250 &16	 &  \textbf{5.54e-05}	 & 6.67e-05	          &     250	 &15	 &  \textbf{8.01e-09}	 & 3.11e-07	                                                        \\ \hline
		300 &150	 &  \textbf{2.11e-10}	 & \textbf{2.11e-10}  &    300	 &28	 &  \textbf{6.11e-10}	 & 1.94e-07	                                                          \\ \hline
	\end{tabular}}
	\caption{Example 5, $t\in \{0.1,\ 0.2,\ 0.5,\ 0.8 \}.$ }
\end{table}

In order to justify this behaviour, we have plotted the error curves $\textrm{err}_{L}(\tau):=f(\tau)-L_n(w, f,\tau)$, and $\textrm{err}_{VP}(\tau):=f(\tau)-V_m^n (w,f,\tau)$ with $\tau\in [-1,1],$  in the cases $n=200,\ m=35$ (Fig. 2, left) and  $n=250,\ m=16$ (Fig. 3, left). \newline
Note that the  selected $n,m$ are those shown in Table 5 for the quadrature errors of the Hilbert transform at $t=0.2$ and $t=0.5$ (see the third and second last lines of the corresponding sub--tables, respectively).
\newline
On the right of each figure, we have zoomed the same plot given on the left in a range close to $\tau=0.2$ (Fig. 2)  and around $\tau=0.5$ (Fig. 3).

As the graphics reveal, the trend of the quadrature errors for $t=0.2$ and $t=0.5$ reflect  that one of the pointwise Lagrange and de la Vall\'ee Poussin approximation errors at $\tau=0.2$ and $\tau=0.5$, respectively, in agreement with the theoretical error estimates (\ref{errore_quad2}) and (\ref{errore_quad3}).

\begin{figure}[h]
	\begin{center}
		\includegraphics[scale=.30]{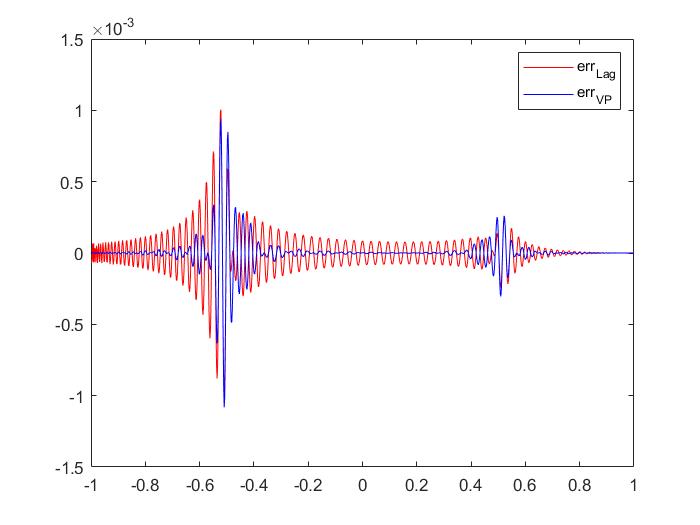}
		\includegraphics[scale=.30]{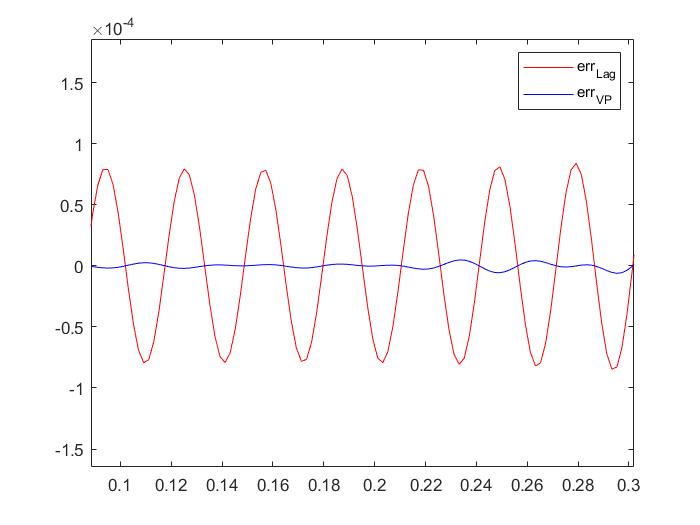}
		{\caption{\label{funzione1}Example 5: graphics of the errors $\textrm{err}_{L}(\tau)$  and $\textrm{err}_{VP}(\tau)$  for $n=200,\ m=35$ (on the left) and zooming of the same plot around $\tau=0.2$ (on the right)}}
	\end{center}
\end{figure}

\begin{figure}[h]
	\begin{center}
		\includegraphics[scale=.30]{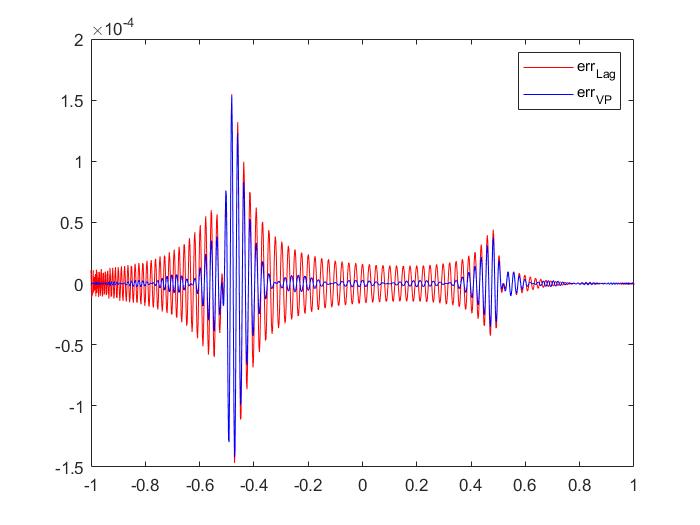}
		\includegraphics[scale=.30]{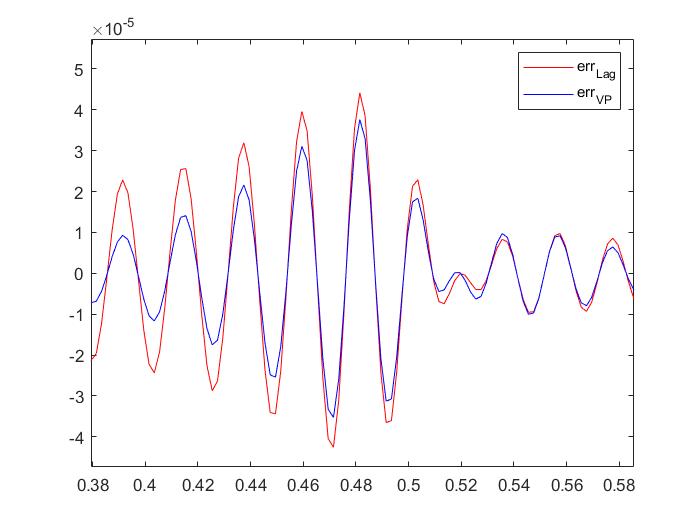}
		{\caption{\label{funzione2}Example 5: graphics of the errors $f(\tau)-L_n (w,f,\tau)$  and $f(\tau)-V_m^n (w,f,\tau)$ for $n=250,\ m=16$ (on the left) and zooming of the same plot around $\tau=0.5$ (on the right)}}
	\end{center}
\end{figure}

\vspace{0.5cm}

\noindent{\bf Test A}

In this test,  we compare the errors obtained by several VP-rules $\{\mathcal{H}^u_{n,m}(w,f,t)\}_{n}$, that differ for the choice of the additional parameter $m$. More precisely, we choose $m=\lfloor n\theta\rfloor$ and let $\theta$ varying in  $(0,1)$, taking also the limit case $\theta=0$, corresponding  to the  L-rule.
On $\theta$ varying in the open interval $(0,1)$,
we have  different VP-rules sequences, all of them convergent, but allowing to  different errors.

Consider indeed the following example
\[\mathcal{H}^u(f,t)=\int_{-1}^1\frac{|x|^\frac 1 7}{x-t }\sqrt[4]{\frac{1-x}{1+x}}dx. \]
$$ u=v^{\frac 1 4 ,-\frac 1 4 }, \ w=u, \ u_1=v^{\frac 1 4,0}\ \ f\in B_r(u_1),\quad r<\frac 1 7.$$
In this case, the function is very smooth, except that in a small interval around  $0$ where a cusp holds.
Here the exact value has been computed by the MG rule with $n=1200$.
In  Figure 4 we show for a fixed $t$, the graphics of the errors in log-scale, for increasing values of $n$ and for fixed   ratios $\theta=\frac m n \in \{0.3, \ 0.5,\ 0.7, 0.9\}.$

As we can see,   all the sequences $\{\mathcal{H}^u_{n,m}(w,f,t)\}_{n}$ converge faster than the  L-rule, for any choice of $\theta$.

\begin{figure}[h]
	\begin{center}
		\includegraphics[scale=.30]{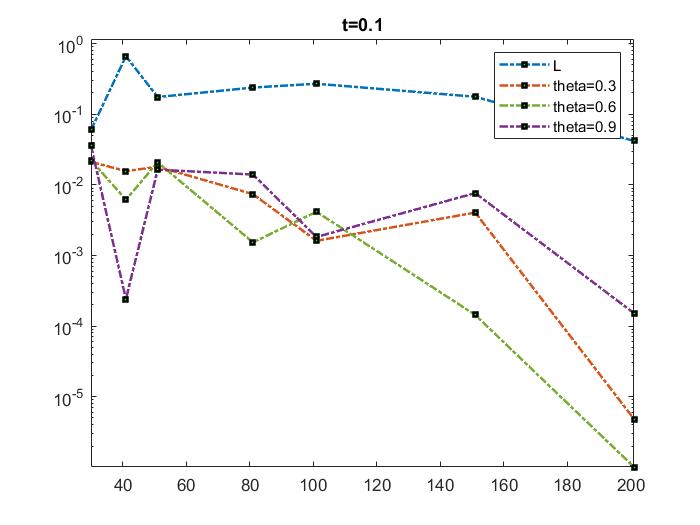}
		\includegraphics[scale=.30]{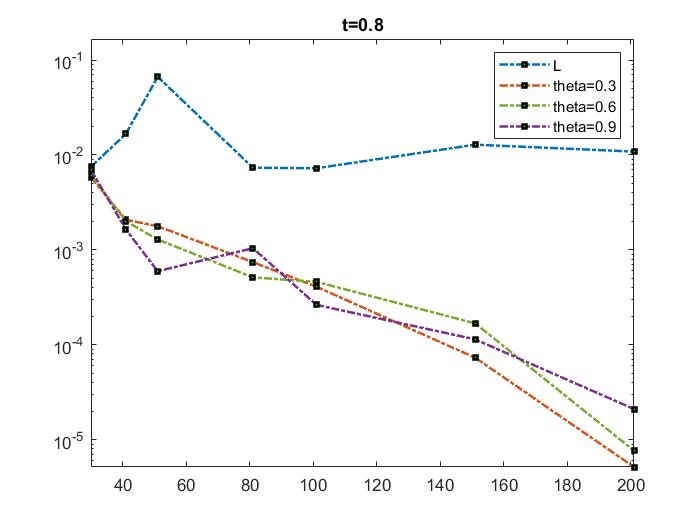}
		{\caption{\label{errori_test7}Absolute errors for Test A, $t=0.1$ (left), $t=0.8$ (right)}}
	\end{center}
\end{figure}

{\bf Test B}

As we have remarked  in the previous test, the rule $\mathcal{H}^u_{m,n}(f)$ for $n$ fixed  and $m$ varying, induces different errors. In this experiment, for any fixed $n$, we determine
the value of $m$ for which the smallest error  is attained  in  VP-rule (say it the ``optimal'' $m$) and for the same $n$, the value of $m$ for which the largest error is achieved (say it the ``worst''  $m$).
Referred to the previous  Example 5, we produce for increasing values of $n$ (in linear scale)  the plot of  the absolute errors  (in log-scale) of VP-rules (either  the best and the worst) and those  due to the  L-rule and the MG-rule.
Here, as well as in Example 5,   we have retained exact the values achieved  by means of the L-rule of order $n=1000.$

\begin{figure}[h]
	\begin{center}
		\includegraphics[scale=.30]{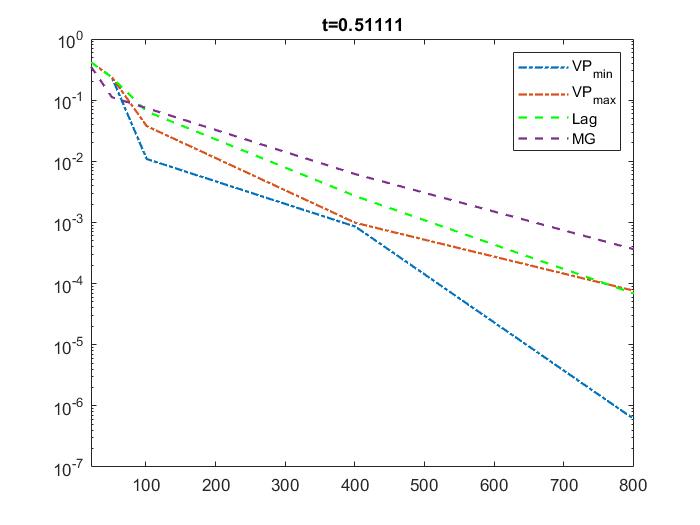}
		\includegraphics[scale=.30]{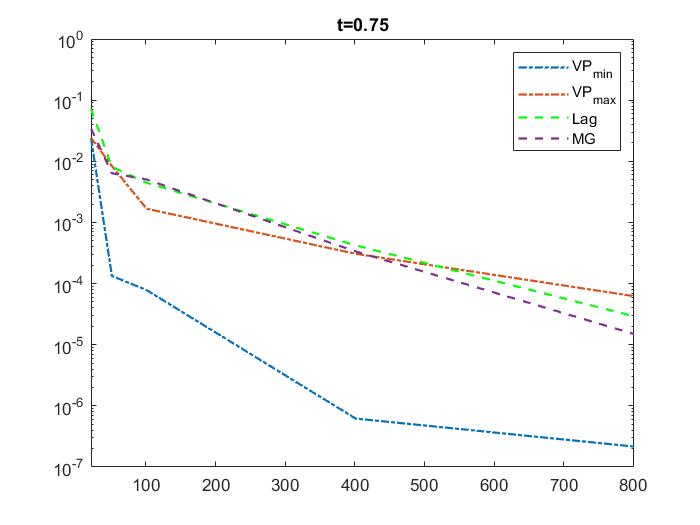}
		{\caption{\label{funzione3}Absolute errors for Example 5, $t=0.51111$ (left), $t=0.75$ (right)}}
	\end{center}
\end{figure}

The errors behaviour confirm that the results by the VP-rule are globally better than those achieved by  the L-rule, when the function $f$ presents a quick variation in  a localized range. As matter of fact, the parameter  $m$ can be used to reduce the quadrature error.

\section{Proofs}
{\bf Proof of Lemma \ref{lem-u}}

In the case $r=0$, equivalence (\ref{equi-lemma})  was proved in \cite[Eq.(2.8)]{Uwe2005}. The proof was mainly based on the following crucial  equivalences
\begin{equation*}
	\sum_{n=1}^\infty \frac{E(n)}n\ \sim\ \sum_{j=0}^\infty E(2^j)\ \sim \ E(1)+\int_0^1\frac{E(t^{-\theta})}tdt,
\end{equation*}
where $E(x)$ stands for any nonnegative decreasing function defined on $[1,\infty)$ and  $\theta$ is an arbitrary positive fixed number.

Therefore, in the general  case $r>0$, (\ref{equi-lemma})  can be proved with analogous arguments by replacing the previous equivalences with the following ones
\begin{equation*}
	\sum_{n=1}^\infty n^{r-1}E(n)\ \sim\ \sum_{j=0}^\infty 2^{jr}E(2^j)\ \sim \ E(1)+\int_0^1\frac{E(t^{-\theta})}{t^{r+1}}dt.
\end{equation*}
\Proofend
{\bf Proof of Lemma \ref{lemmaEm}}

The statement easily follows from the Jackson type inequalities in (\ref{direct}). Using the first one of such inequalities and recalling that $\omega_\varphi^k(f,t)_v$ is an increasing function of $t$, for all $f\in B_0(v)$ we get
\begin{eqnarray*}
E_n(f)_v &\leq& \C\omega^k_{\varphi}\left(f,\frac 1 n\right)_v=
\frac{\C}{\log n}\omega^k_{\varphi}\left(f,\frac 1 n\right)_v\int_\frac 1n^1\frac{dt}t\\
&\le&
\frac{\C}{\log n} \int_{\frac 1 n}^1 \frac{\omega^k_{\varphi}\left(f,t\right)_v}{t}\ dt \leq \frac{\C}{\log n} \|f\|_{B_{0}(v)},
\end{eqnarray*}
having used  (\ref{normB0-equi}) in the last step.

Moreover, starting from the second inequality in (\ref{direct}), for any $f\in B_{r}(v)$ and $k>r$, we get
\[
E_n(f)_v \leq \C \int_0^\frac{1}{n}\frac{\Omega_\varphi^k(f,t)_v}{t}\ dt=
 \C \int_0^\frac{1}{n}\frac{t^r\Omega_\varphi^k(f,t)_v}{t^{r+1}}\ dt
\leq \frac{\C}{n^r} \int_0^\frac{1}{n}\frac{\Omega_\varphi^k(f,t)_v}{t^{r+1}}\ dt
\]
and the second inequality in (\ref{stimaEm}) follows by (\ref{equi-B}).
\Proofend
\textbf{Proof of Theorem \ref{th-Hilbert}}

We remark that in \cite[Th. 3.1]{LutherRusso} the authors have already proved that
 \[
 | \H^uf(t)|u_2(t)\le \C \|f\|_{B_0(u_1)}, \qquad \C\ne \C(f,t)
 \]
 holds. Here we show how the slightly different  estimate in (\ref{stimaUR}) can easily deduced from the proof given in \cite{LutherRusso}.

 Indeed, starting from the following decomposition
 \begin{eqnarray*}
	\mathcal{H}^uf(t)&=& \int_{-1}^1 \frac{f(x)}{x-t} \frac{u_+(x)}{u_-(x)} dx\\
&=&f(t)u_+(t) \int_{-1}^1 \frac{u_-^{-1}(x)}{x-t} \ dx
	+ \int_{-1}^1\frac{f(x)u_+(x)-f(t)u_+(t)}{x-t}\frac{dx}{u_-(x)} \\
	&=:& J_1(t)+J_2(t),
\end{eqnarray*}
the first term $J_1$ can be estimated by means of \cite[Lemma 4.4]{LutherRusso}, which yields
\[
\left|\int_{-1}^1 \frac{u_-^{-1}(x)}{x-t} \ dx\right|\leq \C u_2^{-1}(t), \quad \mathcal{C}\neq{\mathcal{C}(f,t)},
\]
and consequently
\[
|J_1(t)| u_2(t)\le \C |f(t)| u_1(t), \qquad \C\ne \C(f,t).
\]
The estimate of the other term $J_2$ can be achieved taking into account that
\begin{eqnarray*}
	J_2(t)&=& \left( \int_{-1}^{\frac{t-1}{2}}+\int_{\frac{t+1}{2}}^{1}\right)\frac{f(x)u_+(x)-f(t)u_+(t)}{x-t}u_-^{-1}(x)\ dx\\
	&+& \int_{\frac{t-1}{2}}^{\frac{t+1}{2}}\frac{f(x)u_+(x)-f(t)u_+(t)}{x-t}u_-^{-1}(x)\ dx \\
	&=:& I_2+I_3
\end{eqnarray*}
and applying the estimates of these two addenda given in the proof of
\cite[Th. 3.1]{LutherRusso}.
\Proofend
{\bf Proof of Theorem \ref{th-VP_Besov}}

Let us first prove that for all $r\ge 0$, the map $V_n^m(w):B_r(v)\rightarrow B_r(v)$ is uniformly bounded w.r.t. $n$. To this aim we observe that
\[
E_k(V_n^m(w,f))_v=0, \qquad \forall k\ge n+m-1,
\]
since $V_n^m(w,f)\in \PP_{n+m-1}$. Moreover, for $k=0,\ldots, n+m-2$, we note that
\[
E_k(V_n^m(w,f))_v\le E_k(V_n^m(w,f)- f)_v+ E_k(f)_v\le \|[f-V_n^m(w,f)]v\|_\infty +  E_k(f)_v.
\]
Consequently, for all $r\ge 0$, since we are assuming that $V_n^m(w)$ is uniformly bounded in $C^0_v$, we get
\begin{eqnarray*}
& &\|V_n^m(w,f)\|_{B_r(v)}=\|V_n^m(w,f)v\|_\infty+\sum_{k=0}^{n+m-2}  (k+1)^{r-1}E_k(V_n^m(w,f))_v\\
&\le& \C \|fu\|_\infty + \|[f-V_n^m(w,f)]v\|_\infty\sum_{k=0}^{n+m-2}  (k+1)^{r-1}+\sum_{k=0}^{n+m-2}  (k+1)^{r-1}E_k(f)_v\\
&\le& \C \|f\|_{B_r(v)} + \|[f-V_n^m(w,f)]v\|_\infty\sum_{k=0}^{n+m-2}  (k+1)^{r-1}
\end{eqnarray*}
Hence, taking into account that
\[
\sum_{k=0}^{n+m-2}  (k+1)^{r-1}\le \C \left\{\begin{array}{ll}
\log n, & \mbox{if $r=0$}\\
n^r, & \mbox{if $r>0$}
\end{array}\right., \qquad \C\ne\C(n),
\]
from (\ref{stimaVP_1})  we deduce
\begin{equation*}\label{thesi}
\|V_n^m(w,f)\|_{B_r(v)}\le \C \|f\|_{B_r(v)} ,\quad \C\ne\C(n,f), \qquad \forall r\ge 0.
\end{equation*}
This equation and (\ref{inva}), for all $P\in\PP_{n-m}$, yield
\[
\|f-V_n^m(w,f)\|_{B_r(v)}\le \|f-P\|_{B_r(v)}+\|V_n^m(w,f-P)\|_{B_r(v)}\le \C\|f-P\|_{B_r(v)},
\]
where $\C\ne\C(n,f,P)$ and (\ref{eq-conv}) follows by taking into account that $(n-m)\sim n$ and that the polynomials are dense in the space $B_r(u)$, $r\ge 0$ (see \cite{Uwe2005}, \cite{MR_Besov_inf}).

Finally, let us prove  (\ref{VP_map1})--(\ref{VP_map2}).

We remark that
\begin{equation*}\label{erroreBA}
E_k(f-V_n^m(w,f))_v  \left[\begin{array}{lr}
 =E_k(f)_v,  &  \mbox{if \ $k\geq n+m-1$}, \\
 & \\
\leq  \|[f-V_n^m(w,f)]v\|_\infty, & \mbox{if \ $k<n+m-1$}.
\end{array}
\right.
\end{equation*}
Consequently, for all $s\ge 0$, we have
\begin{eqnarray*}
& &\|f-V_n^m (w, f)\|_{B_s(v)}  = \|[f-V_n^m(w,f)]v\|_\infty
 + \sum_{k=0}^\infty (k+1)^{s-1}E_k(f-V_n^mf)_v\\
 & \leq& \|[f-V_n^m(w,f)]v\|_\infty
\left[1+\sum_{k=1}^{n+m-2}(k+1)^{s-1}\right]
 + \sum_{k=n+m-1}^\infty (k+1)^{s-1}E_k(f)_v
\end{eqnarray*}
Therefore, using (\ref{stimaVP_1}) and (\ref{stimaEm}), in the case $s=0$ we get
\begin{eqnarray*}
\|f-V_n^m (w, f)\|_{B_0(v)}&\le & \C \frac{\|f\|_{B_r(v)}}{n^r}\left[1+\sum_{k=1}^{n+m-2}\frac 1{k+1}\right]\\
&+&\C \|f\|_{B_r(v)} \sum_{k=n+m-1}^\infty \frac 1{(k+1)^{r+1}}\le  \frac{\|f\|_{B_r(v)}}{n^r} (\C+\C\log n)
\end{eqnarray*}
as well as, for all $s>0$, we have
\begin{eqnarray*}
\|f-V_n^m (w, f)\|_{B_s(v)}&\le& \C \frac{\|f\|_{B_r(v)}}{n^r}\left[1+\sum_{k=1}^{n+m-2} (k+1)^{s-1}\right]\\
&+ &\C \|f\|_{B_r(v)} \sum_{k=n+m-1}^\infty (k+1)^{s-r-1} \leq \frac{\|f\|_{B_r(v)}}{n^{r}} \left(\C+ \C n^s\right)
\end{eqnarray*}
\Proofend
{\bf Proof of Theorem \ref{convLag}}

We point out that the equivalence between (\ref{condConvLag}) and (\ref{lagrange}) can be found in \cite[Th. 4.3.1]{mastromilobook} while (\ref{Lag_map2}) has been stated in \cite{MR_Besov_inf}. Finally, the proof of (\ref{Lag_map1}) can be carried out by similar arguments used in the proof of Theorem \ref{th-VP_Besov} about $V_n^m(w,f)$. Indeed by (\ref{lagrange}) and (\ref{stimaEm}) we deduce
\begin{eqnarray*}
	\|f-L_n(w, f)\|_{B_0(v)}&\le & \C \|[f-L_n(w,f)]v\|_\infty\left[1+\sum_{k=1}^{n-1}(k+1)^{-1}\right] + \sum_{k=n}^\infty \frac{E_k(f)_v}{k+1}\\
		& \leq & \C \|f\|_{B_r(v)} \frac{\log n}{n^r} \left[1+\sum_{k=1}^{n-1}(k+1)^{-1}\right] + \C \frac{\|f\|_{B_r(v)}}{n^r}\\
		&\leq &  \C \frac{\log^2 n}{n^r}\|f\|_{B_r(v)}.
	\end{eqnarray*}
\Proofend
{\bf Proof of Proposition \ref{th-recQ}}

We start recalling that for $u=v^{a,b}$ it results (\cite[p.320]{Gradshteyn}):
\begin{eqnarray*}\mathcal{H}^u \mathbf{1}(t)&=&\int_{-1}^1\frac{u(x)}{x-t}dx=u(t)\pi\cot \pi(a+1)\\&-&\frac{2^{a+b}\Gamma(a)\Gamma(b+1)}{\Gamma(a+b+1)} \ _2F_1\left(-a-b,1;1-a;\frac{ 1-t} 2\right),
	\end{eqnarray*}
where $_2F_1$ is the generalized hypergeometric function of order $2,\ 1.$
\vspace{0.3cm}
By the recurrence relation (\ref{3term-pol}), it easily follows
\begin{eqnarray*}
      Q_0(t)&=&\frac {\mathcal{H}^u \mathbf{1}(t)}{\sqrt{\int_{-1}^1 u(x)dx}}\\   b_1 Q_1(t)&=& (t-a_0)Q_0(t)+d_0 \\
        b_{j+1}Q_{j+1}(t)&=&(t-a_{j})Q_{j}(t)- b_j Q_{j-1}(t)+d_j,\ \ j\ge 1,
    \end{eqnarray*}
    where for $j=0,1,\dots, $ $d_j=\int_{-1}^1 p_j(w,x) u(x)dx.$
Therefore, setting
$$A_0=b_1^{-1},\ B_0=\frac{a_0}{b_{1}},\ \  $$
$$A_j=b_{j+1}^{-1},\ B_j=-\frac{a_j}{b_{j+1}},\ C_j=\frac{b_j}{b_{j+1}}\ \
\ D_j=\frac{d_j}{b_{j+1}},\quad j=0,1,\dots$$
we have
\begin{eqnarray*}
      Q_0(t)&=&\frac {\mathcal{H}^u \mathbf{1}(t)}{\sqrt{\int_{-1}^1 u(x)dx}},\quad   Q_1(t)= (A_0t+B_0)Q_0(t)+D_0 \\
        Q_{j+1}(t)&=&(A_jt+B_{j})Q_{j}(t)- C_j Q_{j-1}(t)+D_j,\ \ j\ge 0.
    \end{eqnarray*}
\Proofend
{\bf Proof of Theorem \ref{theo2}.}
 We remark that the definition (\ref{def-err}), (\ref{def-gen}) and (\ref{iniziale}) yield
\[
\E_{n,m}^u(w,f,t):=\left|\H^uf(t)-\mathcal{H}_{n,m}^u(w,f,t)\right|=
\left|\H^u[f-V_n^m(w,f)](t)\right|, \quad t\in (-1,1).
\]
Therefore, by applying  Theorem \ref{th-Hilbert}, for all $t\in (-1,1)$, we immediately get (\ref{errore_quad_point}).
Consequently (\ref{conv-B0}) follows from (\ref{nearbest}) and (\ref{eq-conv}).
\Proofend
{\bf Acknowledgments}
This research has been accomplished within Rete ITaliana di Approssimazione (RITA) and TA-UMI. The first two authors are partially supported by University of Basilicata (local research funds) and the GNCS-INdAM funds 2020, project ``Approssimazione multivariata ed equazioni funzionali per la modellistica numerica''.

{\bf References}
\bibliographystyle{plain}
\bibliography{biblioHilb}
\end{document}